\documentclass[12pt,reqno]{amsart}

\theoremstyle{plain}
\newtheorem{thm}{Theorem}[section]
\newtheorem{lem}{Lemma}[section]
\newtheorem{prop}{Proposition}[section]

\theoremstyle{definition}

\newcommand{\NN}{\mathbb{N}}

\newcommand{\RR}{\mathbb{R}}

\newcommand{\bu}{\mathbf{u}}
\newcommand{\bv}{\mathbf{v}}

\newcommand{\bx}{\mathbf{x}}

\newcommand{\bbf}{\mathbf{f}}

\newcommand{\bnabla}{\boldsymbol{\nabla}}

\newcommand{\define}{\stackrel{\text{\rm def}}{=}}
\newcommand{\Ccal}{{\mathcal C}}
\newcommand{\Acal}{{\mathcal A}}

\newcommand{\Gcal}{{\mathcal G}}

\newcommand{\Vcal}{{\mathcal V}}

\newcommand{\reg}{{\text{\rm reg}}}

\newcommand{\Ocal}{{\mathcal O}}

\newcommand{\rd}{{\text{\rm d}}}
\newcommand{\rw}{{\text{\rm w}}}

\newcommand{\ddt}[1]{\frac{\text{\rm d}#1}{\text{\rm dt}}}

\newcommand{\Vinner}[1]{(\!({#1})\!)}

\begin{document}
\numberwithin{equation}{section}

\title[Asymptotic conditions for strong convergence to the 3D NSE 
  attractor]{Asymptotic regularity
  conditions for the strong convergence towards weak limit sets and 
  weak attractors of the 3D Navier-Stokes equations}
\date{October 27, 2003}

\author{Ricardo M. S. Rosa}
\address[Ricardo M. S. Rosa]{Departamento de Matem\'atica Aplicada \\ 
  Instituto de Matem\'atica \\ Universidade Federal do Rio de Janeiro \\ 
  Caixa Postal 68530 Ilha do Fund\~ao \\ Rio de Janeiro RJ 21945-970 \\
  BRAZIL}
\email{rrosa@ufrj.br}
\thanks{The author was partly supported by a fellowship from 
  CNPq, Bras\'{\i}lia, Brazil}


\begin{abstract}
  The asymptotic behavior of solutions of the three-dimensional 
  Navier-Stokes equations is considered on bounded smooth domains with
  no-slip boundary conditions or on periodic domains.
  Asymptotic regularity conditions are presented to ensure that the 
  convergence of a Leray-Hopf weak solution to its weak $\omega$-limit 
  set -- weak in the sense of the weak topology of the space $H$ of 
  square-integrable divergence-free velocity fields -- are achieved also 
  in the strong topology of $H$. In particular, if a weak $\omega$-limit 
  set is bounded in the space $V$ of velocity fields 
  with square-integrable vorticity then the attraction to the set holds 
  also in the strong topology of $H$. Corresponding results for the 
  strong convergence towards the weak global attractor of Foias and Temam
  are also presented.
\end{abstract}

\maketitle

\section{Introduction}  

The notions of limit sets and attractors (whether local or global) permeate 
the theory of dynamical systems both in finite and in infinite dimensions. 
In the case of infinite dimensions, the existence of such sets, 
in particular that of the global attractor, is a major issue 
(the global attractor is the minimal set for the 
inclusion relation which uniformly attracts all bounded sets of initial 
conditions -- the global attractor contains all locally attracting sets and 
$\omega$-limit sets).
Existence results have been obtained for a number of nonlinear 
partial differential equations modeling various phenomena. 

In this note we address the celebrated system associated with the 
Navier-Stokes equations for an incompressible fluid filling a region 
in a three-dimensional space (3D NSE for short). 
Due in particular to the lack of a result on the global well-posedness for 
the 3D NSE the notion of attractor in this case is not settled, and the
study of the asymptotic behavior of this system is a major challenge.

In \cite{foiastemam} Foias and Temam introduced a notion of {\em 
weak} global attractor (see the definition in \eqref{weakglobalattractor}), 
which is loosely speaking a global attractor for
the weak topology of the natural phase space of square-integrable 
divergence-free vector fields (see also the related notion of trajectory 
attractor \cite{ball2,chepyzhovvishik,sell}).
A notion of weak limit set (limit set for the weak topology, see definition
\eqref{weakomegalimitset})) can be 
similarly considered. The study of attractors and limit sets for the 
strong topology is more delicate due to the lack of global regularity
and uniqueness of the solutions.

Our aim in this note is to consider weak limit sets
and the weak global attractor of Foias and Temam 
and present a condition for convergence
in the strong topology of the phase space. It is an asymptotic regularity 
condition.
More precisely, we prove that if the weak $\omega$-limit set (respectively 
weak global attractor) is made of points through which pass solutions 
satisfying the energy {\em equation} (with equality, not just inequality, 
see \eqref{energyequation}), then it is a strong $\omega$-limit
set (resp. global attractor for the strong topology). 
A corollary of this results with a simpler condition for
strong convergence is that the weak $\omega$-limit set
(resp. weak global attractor) be bounded in the space 
of divergence-free velocity fields with finite enstrophy (i.e. 
square-integrable vorticity). 
This includes weakly attracting fixed points, weak limit cycles, 
weakly attracting quasi-periodic orbits, etc,
which turn out to be strongly attracting provided they have bounded
enstrophy. It also includes hyperbolic objects,
in which case the weakly attracting invariant manifolds turn
out to be also strongly attracting provided they have bounded enstrophy.

The proof of this result is based on an idea devised by Ball 
(see \cite{ball,ball3,ghidaglia2,rosa1,mrw}),
exploting energy-type equations to prove asymptotic compactness of the
trajectories. Let us consider limit sets for simplicity. We start
with a trajectory $\bu=\bu(t)$. This trajectory is said to be
asymptotically compact (in a given space) if given any time 
sequence $t_n \rightarrow \infty$, there exists a convergent subsequence 
for $\{\bu(t_n)\}_n$. This asymptotic compactness implies the existence
of the corresponding $\omega$-limit set. A similar result holds for
global attractors \cite{lady, hale, temam2, sellyou}.

The idea of the energy-equation
method to obtain the asymptotic compactness can be broken down into
two steps: i) weak compactness of the sequence $\{\bu(t_n)\}_n$, and
ii) norm convergence $|\bu(t_{n_j})| \rightarrow |\bv_0|$ 
of a weakly convergente subsequence $\bu(t_{n_j})\rightharpoonup \bv_0$, 
as $j\rightarrow \infty$. In uniformly convex spaces (such as Hilbert
spaces), weak plus norm convergences implies strong convergence, hence
asymptotic compactness in the strong topology. In pratice, the
first step follows from classical a~priori estimates obtained from energy-type 
{\em inequalities}, while the second one, as developed in \cite{ball,ball3}, 
follows from energy-type {\em equations} (the equality is important!).

At this point one may be skeptical about the condition of an energy equation
for the 3D NSE since the (Leray-Hopf) weak solutions are known to satisfy 
only an energy inequality. The crucial point that we present here is based 
on the simple observation that in fact the energy equation need only be 
satisfied by the limit solution associated with $\bv_0$. Hence, asymptotic
regularity is all what is needed for the strong convergence of a 
Leray-Hopf weak solution towards its limit set. The precise results 
are given in Lemma \ref{lemmaforoneelementofomega}, Proposition
\ref{propforomega}, and Theorems \ref{theoremforomegainVreg}
to \ref{theoremforglobalstrongsolutions} below.

It should be clear that this idea may be adapted to yield similar results 
for other differential equations in which uniqueness and lack of regularity 
are troublesome, such as wave equations with critical nonlinearities.

\section{Preliminaries}

We recall now some classical results which can be found for instance
in \cite{lady63, lions, constantinfoias,temam}.
We consider the three-dimensional Navier-Stokes equations with either
periodic or no-slip boundary conditions. In the periodic case, we consider
the whole space $\RR^3$, and the flow is assumed periodic with period $L_i$ 
in each direction $x_i$. We define $\Omega=\Pi_{i=1}^3(0,L_i)$ and assume
that the average flow on $\Omega$ vanishes, i.e.
\[ \int_\Omega \bu(\bx) \;\rd\bx = 0.
\]
Here, $\bu=(u_1,u_2,u_3)$ denotes the velocity vector,
and $\bx=(x_1,x_2,x_3)$, the space-variable.

In the no-slip case, the flow is considered in a bounded domain $\Omega$
of $\RR^3$, with smooth boundary $\partial \Omega$, and the
no-slip boundary condition on $\partial\Omega$ is assumed, 
i.e. $\bu=0$ on $\partial\Omega$.

In either periodic or no-slip case, we obtain a functional equation formulation
for the time-dependent velocity field $\bu=\bu(t)$ in a suitable space $H$:
\begin{equation}
  \label{nseeq}
  \ddt{\bu} + \nu A\bu + B(\bu,\bu) = \bbf.
\end{equation}
We consider the test spaces 
\[ \Vcal=\left\{\bu\in \Ccal^\infty(\Omega); \;\bnabla\cdot\bu=0, 
      \;\bu(\bx) \text{ is periodic with period } L_i \text{ in each 
       direction $x_i$} \right\}, 
\]
in the periodic case, and 
\[ \Vcal=\left\{\bu\in \Ccal_\textrm{c}^\infty(\Omega)^3; \;
      \bnabla\cdot\bu = 0 \right\}, 
\]
in the no-slip case, 
and define $H$ as the completion of $\Vcal$ under the $L^2(\Omega)^3$ norm
(where $\Ccal_\textrm{c}^\infty(\Omega)$ denotes the space of 
infinitely-differentiable real-valued functions with compact support 
on $\Omega$). 
We also consider the space $V$ defined as the completion of $\Vcal$ under
the $H^1(\Omega)^3$ norm. We identify $H$ with its dual and consider
the dual space $V'$, so that $V\subset H \subset V'$. We denote by
$H_\rw$ the space $H$ endowed with its weak topology.

We consider the inner products in $H$ and $V$ given respectively by
\[ (\bu,\bv) = \int_\Omega \bu(\bx)\cdot\bv(\bx) \;\rd\bx,
     \quad \Vinner{\bu,\bv} = \int_\Omega \sum_{i=1,3}
        \frac{\partial \bu}{\partial x_i} 
          \cdot \frac{\partial \bv}{\partial x_i}\; \rd\bx,
\]
and the associated norms $|\bu|=(\bu,\bu)^{1/2}$, 
$\|\bu\|=\Vinner{\bu,\bu}^{1/2}$. 

We denote by $P_{\text{LH}}$ the (Leray-Helmhotz) orthogonal 
projector in $L^2(\Omega)^3$ onto the subspace $H$. In
\eqref{nseeq}, $A$ is the Stokes operator $A\bu = - P_{\text{LH}}\Delta \bu$;
$B(\bu,\bv) = P_{\text{LH}}((\bu\cdot\bnabla)\bv)$ is a bilinear term 
corresponding to the inertial term; $\bbf$ 
represents the mass density of volume forces applied to the fluid, 
and we assume that $\bbf\in H$; and $\nu>0$ is the
kinematic viscosity. 
The Stokes operator is a positive self-adjoint operator on $H$, 
and we denote by $\lambda_1>0$ its first eigenvalue.

A {\em Leray-Hopf weak solution} on an open time interval $I=(t_0,t_1)$,
$-\infty\leq t_0 < t_1 \leq \infty$,
is defined as a function $\bu=\bu(t)$ on $(t_0,t_1)$ with values in $H$
and satisfying the following properties:
\begin{itemize}
  \item[(i)] $\bu\in L^\infty(t_0,t_1;H)\cap L_{\text{loc}}^2(t_0,t_1;V)$;
  \item[(ii)] $\partial\bu/\partial t \in L_{\text{loc}}^{4/3}(t_0,t_1;V')$;
  \item[(iii)] $\bu\in \Ccal(I;H_\rw)$ (i.e. weakly continuous);
  \item[(iv)] $\bu$ satisfies the functional equation \eqref{nseeq}
    almost everywhere on $I=(t_0,t_1)$;
  \item[(v)] $\bu$ satisfies the following energy inequality in the 
    distribution sense on $I=(t_0,t_1)$:
      \begin{equation}
        \label{energyinequality}
          \frac{1}{2}\frac{\rd}{\rd t}|\bu(t)|^2 + \nu \|\bu(t)\|^2 
            \leq (\bbf,\bu(t)).
      \end{equation}
\end{itemize}

A {\em Leray-Hopf weak solution} on an interval of the form $[t_0,t_1)$ is 
defined as a Leray-Hopf weak solution on $(t_0,t_1)$ which is continuous
at $t=t_0$, i.e. 
\begin{itemize}
  \item[(vi)] $\bu(t)\rightarrow \bu(t_0)$ in $H$, as $t\rightarrow t_0^+$.
\end{itemize}

From now on, for notational simplicity, a weak solution will always
mean a Leray-Hopf weak solution. For a weak solution on an arbitrary
interval $I$, it follows that
\begin{equation}
  \label{energyestimate}
  |\bu(t)|^2 \leq |\bu(t')|^2 e^{-\nu\lambda_1 (t-t')}
     + \frac{1}{\nu^2\lambda_1^2} |\bbf|^2 
         \left(1-e^{-\nu\lambda_1 (t-t')}\right),
\end{equation}
for all $t$ in $I$ and almost all $t'$ in $I$ with $t'<t$. The allowed
times $t'$ are the Lebesgue points of the function $t\mapsto |\bu(t)|$.
In the case of a weak solution on an interval of the form $[t_0,t_1)$,
the point $t_0$ is a point of continuity of $t\mapsto |\bu(t)|$, hence
a Lebesgue point, so that the estimate above is also valid for the initial
time $t'=t_0$. 

Another classical estimate obtained from the energy inequality 
is 
\begin{equation}
  \label{L2enstrophyestimate}
  |\bu(t)|^2 + \nu \int_{t'}^t \|\bu(s)\|^2 \;\rd s
      \leq |\bu(t')|^2 + \frac{1}{\nu\lambda_1} |\bbf|^2 (t-t'),
\end{equation}
for all $t$ in $I$ and almost all $t'$ in $I$ with $t'<t$, with the
set of allowed times $t'$ consisting again as the Lebesgue points
of the function $t\mapsto |\bu(t)|$.

A {\em strong solution} on an arbitrary interval $I$ is defined as
a weak solution on $I$ satisfying 
\begin{itemize}
  \item[(vii)] $\bu\in \Ccal(I;V)$.
\end{itemize}
Any strong solution satisfies the energy equation
\begin{equation}
  \label{energyequation}
  \frac{1}{2}\frac{\rd}{\rd t}|\bu(t)|^2 + \nu \|\bu(t)\|^2 
          = (\bbf,\bu(t))
\end{equation}
in the distribution sense on its interval of definition.

It is well established that given any initial time $t_0$ and any initial 
condition $\bu_0$ in $H$, there exists at least one global 
weak solution on $[t_0,\infty)$ satisfying $\bu(t_0)=\bu_0$. 
It is also known that if $\bu_0$ belongs to $V$, then there exists
a local strong solution, defined on some interval $[t_0,t_1)$,
$t_1>t_0$, with $\bu(t_0)=\bu_0$. Finally, any strong solution is unique
on its interval of definition. The uniqueness in this case is with
respect to the larger class of all weak solution, i.e. any weak solution
with $\bu(t_0)=\bu_0$ must agree with the strong solution on the interval
of definition of the latter. 

\section{Weak limit sets and the weak global attractor}

The weak global attractor, as defined in \cite{foiastemam},
is the set
\begin{equation}
  \label{weakglobalattractor}
  \Acal_\rw = \left\{ \bv_0 \in H; \;\parbox{4in}{there exists at least
      one global weak solution $\bv=\bv(t)$ defined for all $t\in\RR$ 
      which is uniformly bounded in $H$, i.e. 
      $\sup_{t\in \RR}|\bv(t)|<\infty$, and such that $\bv(0)=\bv_0$}\right\}.
\end{equation}
Due to the energy estimate \eqref{energyestimate} and the uniform
bound on the global solutions in the definition of $\Acal_\rw$ it follows
that $\Acal_\rw$ is a bounded set in $H$:
\[ |\bv_0| \leq \frac{1}{\nu\lambda_1}|\bbf|, \quad \forall 
      \bv_0\in \Acal_\rw.
\]
It is proved in \cite{foiastemam} that $\Acal_\rw$ is weakly
compact in $H$ and that it attracts all weak solutions in the
following sense: If $\bu=\bu(t)$ is a weak solution on $[t_0,\infty)$
for some $t_0\in \RR$, then for any neighborhood $\Ocal$ of $\Acal_\rw$
in the weak topology of $H$, there exists a time $T\geq t_0$ such 
that $\bu(t)\in \Ocal$ for all $t\geq T$. Since $H$ is separable
the weak topology is metrizable on bounded sets and the convergence
above can be rewritten is terms of this metric. Finally, a certain
invariance property holds, namely if $\bv_0$ belongs to $\Acal_\rw$ and
$\bv=\bv(t)$, $t\in \RR$, is a uniformly bounded global weak solution 
through $\bv_0$ then
$\bv(t)\in \Acal_\rw$ for all $t\in \RR$. 

Besides the pointwise attraction (attraction of individual weak solutions)
of the weak global attractor, one can show that the attraction is in
fact uniform with respect to uniformly bounded sets of initial condition
(see \cite{fmrt}). More precisely, given $t_0\in \RR$ and $R>0$, then for 
every neighborhood $\Ocal$ of $\Acal_\rw$ in the weak topology of $H$, there 
exists a time $T\geq t_0$ such that $\bu(t)\in \Ocal$ for all $t\geq T$ and
for every weak solution $\bu$ on $[t_0,\infty)$ with
$\sup_{t\geq t_0}|\bu(t)| \leq R$. Since $\Acal_\rw$ is bounded in $H$
and the weak topology of $H$ is metrizable on bounded subsets this 
uniform attraction in the weak topology can be rewritten in terms of 
an associated metric. 

These properties define 
$\Acal_\rw$ and justify its definition as the weak global attractor.
They can also be used to characterize $\Acal_\rw$ in a more classical
way:
\begin{equation}
  \label{Acalwasomegalimit}
  \Acal_\rw = \left\{ \bv_0\in H; 
      \;\parbox{4.2in}{there exist global weak solutions 
      $\bu_n=\bu_n(t)$, $n\in \NN$, defined for all $t\geq 0$,
      with $\sup_{t\in \RR}|\bu_n(t)|\leq |\bbf|/\nu\lambda_1$, 
      and a time sequence $\{t_n\}_n$, $t_n\geq 0$, $t_n\rightarrow \infty$,
      such that $\bu_n(t_n) \rightharpoonup \bv_0$} \right\}.
\end{equation}

Now, given an arbitrary weak solution $\bu=\bu(t)$ on
an interval of the form $(t_0,\infty)$ or $[t_0,\infty)$, 
for some $t_0\in \RR$, we define its weak $\omega$-limit set by
\begin{equation}
  \label{weakomegalimitset}
  \omega_\rw(\bu)= \left\{ \bv_0\in H; \;\exists \{t_n\}_n, t_n > t_0,
     t_n\rightarrow \infty, \;\bu(t_n)\rightharpoonup \bv_0 \right\},
\end{equation}
where $\rightharpoonup$ denotes the weak convergence in $H$. This
set is always nonempty since $\{\bu(t)\}_{t> t_0}$ is bounded in $H$
(thanks to \eqref{energyestimate}), hence weakly precompact.
Since the weak topology is metrizable on bounded subsets of $H$,
the classical characterization $\omega_\rw(\bu)=\cap_{t\geq 0}
\overline{\cup_{t\geq s} \{\bu(t)\}}^\rw$ holds, where
$\overline{\ \cdot\ }^\rw$ denotes the closure in the weak topology.
Hence, $\omega_\rw(\bu)$ is weakly compact. By classical 
dynamical system arguments one can also show that $\omega_\rw(\bu)$ 
attracts $\bu$ in the
sense that for any weakly open set $\Ocal$ containing $\omega_\rw(\bu)$,
there exists a time $T>t_0$ such that $\bu(t)\in \Ocal$ for all $t\geq T$.

As for the invariance property, it is possible to show that for every
$\bv_0$ in $\omega_\rw(\bu)$, there exists a global weak solution $\bv=\bv(t)$,
$t\in \RR$, with $\bv(0)=\bv_0$ and $\bv(t)\in \omega_\rw(\bu)$ for
all $t\in \RR$. This is achieved by passing to the limit the
solutions $\bu(t_n+\cdot)$ over time intervals $[-T,T]$, for arbitrarily
large times $T$. In fact, classical a~priori estimates (derived from
\eqref{energyestimate} and \eqref{L2enstrophyestimate}) yield that 
$\{\bu(t_n+\cdot)\}_n$ is bounded on $L^\infty(-T,T;H)\cap L^2(-T,T;V)$ 
and that $\{\partial\bu/\partial t(t_n+\cdot)\}_n$ is
bounded on $L^{4/3}(-T,T;V')$, which imply precompactness in
$L^2(-T;T,H)$ and $\Ccal([-T,T],H_\rw)$. A diagonalization process
guarantees convergence on any bounded set to a global weak solution
defined on all $\RR$. However, due to the possible lack of uniqueness 
one cannot
guarantee the invariance for every solution passing through $\bv_0$. 

With the above invariance property in mind, given $\bv_0$ in $\omega_\rw(\bu)$
we defined $\Gcal_\bu(\bv_0)$ as the set of all global weak solutions
$\bv=\bv(t)$, $t\in \RR$, with $\bv(0)=\bv_0$, and such that there exists 
a sequence
$\{t_n\}_n$, $t_n\geq t_0$, $t_n\rightarrow \infty$, with the property
that $\bu(t_n+\cdot)$ converges to $\bv$ in $\Ccal([-T,T],H_\rw)$, for
all $T>0$. 
Note that this implies that $\Gcal_\bu(\bv_0)\subset \omega_\rw(\bu)$,
for all $\bv_0$ in $\omega_\rw(\bu)$.

Finally, by classical a~priori estimates (derived from
\eqref{energyestimate} and \eqref{L2enstrophyestimate}) and Aubin's compactness
theorem the convergence of $\bu(t_n+\cdot)$ to $\bv$ holds weakly-star
in $L^\infty(-T,T;H)$, weakly in $L^2(-T,T;V)$, and strongly
in $L^2(-T,T;H)$. Then, the convergence required in the definition of
$\Gcal_\bu(\bv_0)$ guarantees the uniqueness of the limit and the convergence
without passing to further subsequences.

A similar argument for the weak global attractor yields that for every
$\bv_0$ in $\Acal_\rw$ and every sequences $\{\bu_n\}_n$ and $\{t_n\}_n$
as in the characterization \eqref{Acalwasomegalimit},
with $\bu_n(t_n)\rightharpoonup \bv_0$, 
there exists subsequences $\{\bu_{n_j}\}_j$
and $\{t_{n_j}\}_j$ such that $\bu_{n_j}(t_{n_j}+\cdot)$ converges
weakly-star in $L^\infty(-T,T;H)$, weakly in $L^2(-T,T;V)$, strongly
in $L^2(-T,T;H)$ and strongly in $\Ccal(-T,T;H_\rw)$, for all $T>0$, to a 
global weak solution $\bv=\bv(t)$, with $\bv(t)\in \Acal_w$, for all
$t\in \RR$, and with $\bv(0)=\bv_0$.

\section{Asymptotic regularity conditions for the strong convergence towards
weak limit sets and the weak global attractor}

As mentioned in the introduction the required regularity condition
is that the limit solutions satisfy the energy equation exactly.
More precisely, we have the following result.

\begin{lem}
  \label{lemmaforoneelementofomega}
  Let $\bu$ be a weak solution defined on some interval of the
  form $(t_0,\infty)$ or $[t_0,\infty)$, for some $t_0\in \RR$. 
  Let $\bv_0\in \omega_\rw(\bu)$ and let $\{t_n\}_n$ be such that 
  $t_n>t_0$, $t_n\rightarrow \infty$,
  and $\bu(t_n)\rightharpoonup \bv_0$ weakly in $H$. 
  If there exists a global weak solution $\bv=\bv(t)$, $t\in \RR$,
  such that $\bu(t_n+\cdot)$ converges to $\bv$ in
  $\Ccal([-T,T],H_\rw)$, for all $T>0$, and which satisfies the energy equation
  \eqref{energyequation} on $\RR$, then $\bu(t_n)$ converges 
  strongly in $H$ to $\bv_0$. 
\end{lem}

\proof In what follows, given $T>0$ we restrict ourselves to $n$ such that 
$t_n-T>t_0$. 

The weak solution satisfies the energy inequality \eqref{energyinequality}. 
By adding and subtracting the term $\nu\lambda_1|\bu(t)|^2/2$, we arrive at
the inequality
\begin{equation}
  \label{modifiedenergyinequality}
  \frac{1}{2}\frac{\rd}{\rd t}|\bu(t)|^2 
        + \frac{\nu\lambda_1}{2} |\bu(t)|^2 + \nu [\bu(t)]^2
        \leq (\bbf,\bu(t)),
\end{equation}
where
\begin{equation}
  \label{equivalentnorm}
  [\bu]^2 =  \|\bu\|^2 - \frac{\lambda_1}{2}|\bu|^2,
\end{equation}
with $[\cdot]$ being a norm in $V$ equivalent to $\|\cdot\|$.
By ``multiplying'' \eqref{modifiedenergyinequality} by appropriate
nonnegative test functions approximating $e^{\nu\lambda t}$ and 
``integrating'' (e.g. use $\varphi(t)=\theta_m(t)e^{\nu\lambda_1 t}$
in the distribution formulation of the equation, 
where $\theta_m(t)$ is the continuous piecewise linear function joining 
linearly the values 
$\theta_m(-T)=0$, $\theta_m(-T+1/m)=1$, $\theta_m(-1/m)=1$, 
$\theta_m(0)=0$, and let $m\rightarrow \infty$)
one finds at the limit that
\[ |\bu(t_n)|^2 \leq |\bu(t_n-T)|^2 e^{-\nu\lambda_1 T}
      - 2\int_{-T}^0 e^{\nu\lambda_1 s} \left\{\nu[\bu(t_n+s)]^2
         - (\bbf,\bu(t_n+s)) \right\} \;\rd s,
\]
for almost all $T$ such that $t_n-T>t_0$ (more precisely, it holds
for all Lebesgue points $T$ of the function $T\rightarrow |\bu(t_n-T)|$). 

Passing to the limit as $n$ goes to infinity and using that
$\bu(t_n+\cdot)$ converges to $\bv$ weakly in $L^2(-T,T;V)$
we find
\[ \limsup_{n\rightarrow \infty} |\bu(t_n)|^2 
     \leq \limsup_{n\rightarrow \infty} |\bu(t_n-T)|^2 e^{-\nu\lambda_1 T}
      - 2\int_{-T}^0 e^{\nu\lambda_1 s} \left\{\nu[\bv(s)]^2
         - (\bbf,\bv(s)) \right\} \;\rd s.
\]
Note that we have used that $[\cdot]$ is an equivalent norm for $V$ 
so that
\[ \int_{-T}^0 e^{\nu\lambda_1 s} [\bv(s)]^2 \;\rd s
     \leq \liminf_{n\rightarrow \infty} 
        \int_{-T}^0 e^{\nu\lambda_1 s} [\bu(t_n+s)]^2 \;\rd s.
\]

As for the limit solution $\bv$ it satisfies the energy equation,
so that
\[ |\bv(0)|^2 = |\bv(-T)|^2 e^{-\nu\lambda_1 T}
      - 2\int_{-T}^0 e^{\nu\lambda_1 s} \left(\nu[\bv(s)]^2
         - (\bbf,\bv(s)) \right) \;\rd s.
\]
By subtracting this equality from the inequality for the limsup 
and using that both $\bv(-T)$ and $\bu(t_n-T)$ are bounded in $H$
independently of $T$ and $n$ we find
\[ \limsup_{n\rightarrow \infty} |\bu(t_n)|^2 - |\bv(0)|^2
     \leq 2R e^{-\nu\lambda_1 T},
\]
for some bound $R>0$.
The previous passage should make clear why we need the
equality in the energy equation only for the limit solution.
Now, since $T$ is arbitrary (except for a set of measure zero), 
we may let $T\rightarrow \infty$ to find that
\[ \limsup_{n\rightarrow \infty} |\bu(t_n)|^2 - |\bv(0)|^2 \leq 0.
\]
Since $\bv(0)=\bv_0$, we obtain
\[  \limsup_{n\rightarrow \infty} |\bu(t_n)| \leq |\bv_0|.
\]
On the other hand, since $\bu(t_n)$ converges weakly in $H$ to
$\bv_0$ we have $|\bv_0| \leq \liminf_{n\rightarrow \infty} |\bu(t_n)|$.
Thus $\lim_{n\rightarrow \infty} |\bu(t_n)| = |\bv_0|$,
which together with the weak convergence implies the strong convergence
$\bu(t_n)\rightarrow \bv_0$ in $H$. This concludes the proof.
\qed
\medskip

\begin{prop}
  \label{propforomega}  
  Let $\bu$ be a weak solution defined on some interval of the form
  $(t_0,\infty)$ or $[t_0,\infty)$, for some $t_0\in \RR$. 
  If for every $\bv_0$ in $\omega_\rw(\bu)$ all the global weak solutions
  in $\Gcal_\bu(\bv_0)$ satisfy the energy equation \eqref{energyequation}
  then  $\omega_\rw(\bu)$ attracts $\bu$ in the strong topology of $H$.
\end{prop}

\proof If this were not true we would find a time sequence $\{t_n\}_n$,
with $t_n > t_0$, $t_n\rightarrow \infty$, and such that $\{\bu(t_n)\}_n$
does not have any subsequence converging strongly in $H$. Now, since 
$\{\bu(t_n)\}_n$ is bounded (thanks to \eqref{energyestimate})
it has a weakly convergent subsequence $\{\bu(t_{n_j})\}_{n_j}$, with
a weak limit $\bv_0$ belonging to $\omega_\rw(\bu)$. Since any solution
$\bv$ in $\Gcal_\bu(\bv_0)$ satisfies the condition of Lemma
\ref{lemmaforoneelementofomega} it follows that $\{\bu(t_{n_j})\}_{n_j}$
converges strongly to $\bv_0$, which is a contradiction. Thus,
$\omega_\rw(\bu)$ attracts $\bu$ in the strong topology of $H$.
\qed
\medskip

We now define the set
\[ V_\reg = \left\{ \bv_0\in V; \;\exists \delta>0 \text{ and }
     \exists \text{ a strong solution } \bv\in \Ccal((-\delta,\delta),V) 
     \text{ with } \bv(0)=\bv_0 \right\}.
\]

\begin{thm}
  \label{theoremforomegainVreg}
  Let $\bu$ be a weak solution defined on some interval of the form
  $(t_0,\infty)$ or $[t_0,\infty)$, for some $t_0\in \RR$.
  If $\omega_\rw(\bu)\subset V_\reg$, then $\omega_\rw(\bu)$ attracts
  $\bu$ in the strong topology of $H$.
\end{thm}

\proof In view of Proposition \ref{propforomega} it suffices to show
that any weak solution in $\Gcal_\bu(\bv_0)$ satisfies the energy equation
for all $\bv_0$ in $\omega_\rw(\bu)$. Note first that $\Gcal_\bu(\bv_0)
\subset V_\reg$ since it is included in $\omega_\rw(\bu)$, which is assumed
to be included in $V_\reg$. So any weak solution $\bv$ 
in $\Gcal_\bu(\bv_0)$ is included in $V_\reg$, i.e. $\bv(t)\in V_\reg$ for all 
$t\in \RR$. By the uniqueness of strong solutions among the larger class
of weak solutions it follows that $\bv$ is a global strong solution
and belongs to $\Ccal(\RR,V)$. Therefore, $\bv$ satisfies the energy 
equation, and this completes the proof.
\qed
\medskip

\begin{thm}
  \label{theoremforomegaboundedinV}
  Let $\bu$ be a weak solution defined on some interval of the form
  $(t_0,\infty)$ or $[t_0,\infty)$, for some $t_0\in \RR$.
  If $\omega_\rw(\bu)$ is a bounded subset of $V$, 
  then $\omega_\rw(\bu)$ attracts $\bu$ in the strong topology of $H$.
\end{thm}

\proof
  Following the idea in the proof of Theorem \ref{theoremforomegainVreg},
  any global weak solution $\bv$ in 
  $\Gcal_\bu(\bv_0)$, $\bv_0\in \omega_\rw(\bu)$,
  is uniformly bounded in $V$, hence it is a global strong solution
  belonging to $\Ccal(\RR,V)$, which is sufficient to apply Proposition
  \ref{propforomega} and conclude the proof.
\qed
\medskip

We now present a result concerning the weak global attractor. This
result is not a direct consequence of the previous ones simply because
of the uniform attraction required (reflecting the fact that the weak 
global attractor may be larger than the union of weak $\omega$-limit sets).

\begin{thm}
  \label{theoremforAcalw}
  If $\Acal_\rw \subset V_\reg$ or $\Acal_\rw$ is a bounded subset of $V$, 
  then $\Acal_\rw$ attracts every weak
  solution in the strong topology of $H$, and this attraction 
  is uniform with respect to uniformly bounded sets of weak solutions. 
  More precisely, given $t_0\in \RR$ and $R>0$, then for every
  $\varepsilon > 0$, there exists a time $T\geq t_0$ such that 
  $\text{dist}_H(\bu(t),\Acal_\rw) \define
  \sup_{\bv_0\in \Acal_\rw} |\bu(t)-\bv_0| <
  \varepsilon$, for every weak solution $\bu$ on $[t_0,\infty)$ with
  $\sup_{t\geq t_0}|\bu(t)| \leq R$.
\end{thm}

\proof Suppose the result is not true. Then there exists $t_0\in \RR$,
$R>0$, $\varepsilon>0$, a sequence $\bu_n$ of weak solutions on 
$[t_0,\infty)$ with $\sup_{t\geq t_0}|\bu_n(t)| \leq R$, and a time sequence
$\{t_n\}_n$, $t_n\geq t_0$, $t_n\rightarrow \infty$, such that
\begin{equation}
  \label{contradictionforAcalw}
  |\bu_n(t_n)-\bv_0|\geq \varepsilon, \text{ for all } n 
      \text{ and all } \bv_0\in \Acal_\rw.
\end{equation}

Consider the sequence $\bv_n(t)=\bu_n(t_n+t)$, defined for 
$t\geq t_0-t_n$. By the assumption of uniform estimate on 
$\bu_n$ and by classical a~priori estimates for the 3D NSE 
(derived from \eqref{energyestimate} and \eqref{L2enstrophyestimate}) 
the sequence $\{\bv_n\}_n$ 
is bounded in $L^\infty(-T,T;H)\cap L^2(-T,T;V)$,
and $\{\partial\bv_n/\partial t\}_n$ is
bounded in $L^{4/3}(-T,T;V')$, which imply precompactness in
$L^2(-T;T,H)$ and $\Ccal([-T,T],H_\rw)$. By passing to the limit
as in the classical theory of existence of weak solutions of the
3D NSE and by using a diagonalization process we find a subsequence 
converging on any bounded interval $[-T,T]$, $T>0$, to a global weak 
solution $\bv=\bv(t)$, $t\geq 0$. In particular, 
$\bu_{n_j}(t_{n_j})$ converges weakly to some element $\bv_0=\bv(0)$ in $H$.

At the limit, we retain a uniform bound for $\bv$, 
$\sup_{t\geq t_0}|\bv(t)| \leq R$, so that $\bv(t)$ belongs to
$\Acal_\rw$ for all $t\in \RR$, and in particular $\bv_0\in \Acal_\rw$. 
Since $\Acal_\rw$ either belongs
to $V_\reg$ or is a bounded subset of $V$ it follows 
(see the proofs of Theorems \ref{theoremforomegainVreg} and 
\ref{theoremforomegaboundedinV}) that $\bv$ is a global
strong solution, $\bv\in \Ccal(\RR,V)$, hence $\bv$ satisfies the
energy equation on $\RR$.

Now we proceed as in the proof of Lemma \ref{lemmaforoneelementofomega}.
The only difference being that instead of working with $\bu(t_n+\cdot)$
we work with $\bu_n(t_n+\cdot)$ (we still have the inequalities holding
for almost all $T$, as the countable union of zero-measure sets is still
of zero measure). Then, by repeating the energy-equation argument
with $\bu_n(t_n+\cdot)$ we find that 
$\limsup_{n\rightarrow \infty} |\bu_{n_j}(t_{n_j})| \leq |\bv_0|$, 
so that in fact $\bu_{n_j}(t_{n_j})$ converges strongly to $\bv_0$. But 
this contradicts \eqref{contradictionforAcalw} since $\bv_0$ belongs 
to $\Acal_\rw$. This completes the proof.
\qed
\medskip

Our last result concerns just part of the weak global
attractor. It does not reduce to a statement about
weak $\omega$-limit sets included in the weak global attractor since
it allows for asymptotic limits of sequences of weak solutions
instead of a single solution. It applies, for instance, to
regular (in the sense of being a strong solution) connecting orbits 
which are not necessarily $\omega$-limit sets.

\begin{thm}
  \label{theoremforglobalstrongsolutions}
  Suppose $\bv=\bv(t)$, $t\in \RR$, is a global strong solution uniformly
  bounded in $H$ (in other words, a global weak solution included in 
  $\Acal_\rw\cap V_\reg$) and set $\bv_0=\bv(0)$. Let $\bu_n=\bu_n(t)$, 
  $t\geq 0$, and $\{t_n\}_n$ be as in the characterization 
  \eqref{Acalwasomegalimit}, with $\bu_n(t_n)\rightharpoonup \bv_0$.
  Then, $\bu_n(t_n)\rightarrow \bv_0$ strongly in $H$.
\end{thm}

\proof By assumption, $\bv=\bv(t)$ is a strong solution, hence it satisfies
the energy equation \eqref{energyequation}. 
Then, as in the last paragraph of the proof of 
Theorem \ref{theoremforAcalw} we apply the energy-equation method to
the sequence of weak solutions $\{\bu_n(t_n+\cdot)\}_n$ to show 
that $|\bu_n(t_n)|\rightarrow |\bv_0|$. Then, we conclude that
$\bu(t_n)\rightarrow \bv_0$ strongly in $H$. (At first one may need
to pass to further subsequences, but since the weak limit
$\bu_n(t_n)\rightarrow \bv_0$ exists and hence is unique, the strong 
converge in $H$ must hold for the whole sequence.)
\qed

%
%

\end{document}